\documentclass[12pt,leqno]{article}

\usepackage{amsfonts}
\usepackage{amsmath}
\usepackage{amssymb}
\usepackage{graphicx}

\def\eq{equa\-tion}
\def\eqs{equa\-tions}
\def\ct{con\-ti\-nuous}
\def\fwg{fol\-lowing}
\def\nei{neighbor\-hood}

\def\inr{\mbox{ in }}
\def\onr{\mbox{ on }}
\def\forr{\mbox{ for }}

\def\derp#1#2{\displaystyle\frac{\partial#1}{\partial#2}}
\def\derpp#1{\displaystyle\frac{\partial}{\partial#1}}

\def\<{\left<}
\def\>{\right>}
\def\({\left(}
\def\){\right)}
\def\lac{\left\{}
\def\rac{\right\}}

\def\vsp{\vspace*{1,5mm}\\ }

\def\ff{\forall }

\def\ov{\overline}

\def\wt{\widetilde}

\def\1{^{-1}}
\def\3{\subset }
\def\9{{\infty}}

\def\barr{\begin{array}}
\def\earr{\end{array}}
\def\bit{\begin{itemize}}
\def\eit{\end{itemize}}

\def\G{{\Gamma}}
\def\D{\Delta }
\def\ooo{{\Omega}}
\def\pp{{\partial}}
\def\dd{\displaystyle}
\def\vf{{\varphi}}
\def\lbb{{\lambda}}
\def\g{{\gamma}}
\def\a{{\alpha}}

\def\de{{\delta}}
\def\vp{{\varepsilon}}
\def\na{{\nabla}}

\def\cala{{\mathcal{A}}}

\def\calo{\mathcal{O}}
\def\calf{\mathcal{F}}

\def\call{\mathcal{L}}
\def\calw{\mathcal{W}}

\def\calx{\mathcal{X}}

\def\rr{{\mathbb{R}}}
\def\dd{\displaystyle}
\def\Res{{\rm Re}}
\def\Ims{{\rm Im}}

\def\NS{Navier--Stokes}

\def\one{1\!\!\!\;\mathrm{l}}

\newtheorem{theorem}{Theorem}[section]

\newtheorem{lemma}[theorem]{Lemma}

\newtheorem{remark}{Remark}[section]

\def\pf{\noindent{\bf Proof.} }

\begin{document}

\title{Stabilization of Navier--Stokes equations\\ by oblique boundary feedback controllers}
\author{Viorel Barbu\thanks{Supported by CNCSIS project PNID-/2011.}\\
{\small Octav Mayer Institute of Mathematics of Romanian Academy,
Ia\c si, Romania}}

\maketitle

\begin{abstract}
One designs a linear stabilizable boundary feedback controller for
the Navier--Stokes system on a bounded and open domain $\calo\subset\rr^d$,
$d=2,3$, of the form $u=\eta\dd\sum^N_{j=1}\mu_j\<y,\psi_j\>
\(\derp{\phi_j}n\,(x)+\a(x)\vec n(x)\),$ where $\psi_j,\phi_j$ are related to the eigenfunction system for the adjoint Stokes--Oseen system, $\vec n$ is the
normal to $\pp\calo$, $\<\cdot,\cdot\>$ is the scalar product in$(L^2(\calo))^d$ and $\a$ is any \ct\ function with circulation zero on $\pp\calo$.
\end{abstract}

\section{Introduction}

Consider the \NS\ system
\begin{equation}\label{e1.1}
\barr{ll} \derp Yt-\nu\D Y+(Y\cdot\na)Y=\na p+f_e&\inr\
(0,\9)\times\calo,\vsp \na\cdot\ Y=0&\inr\ (0,\9)\times\calo,\vsp
Y=v&\onr\ (0,\9)\times\pp\calo,\vsp Y(0)=y_0&\inr\
\calo,\earr\end{equation} in a bounded open domain
$\calo\subset\rr^d$, $d=2,3$, with a smooth boundary $\pp\calo$
which, for simplicity, we assume simple connected. Here $\nu>0$,
$f_e$ is a given smooth function and $v$ is a boundary input. If
$y_e$ is an equilibrium solution to \eqref{e1.1}, then
\eqref{e1.1} can be, equivalently, written as
\begin{equation}\label{e1.1prim}
\barr{ll} \multicolumn{2}{l}{\derp Yt-\nu\D
Y+(y_e\cdot\na)Y+(Y\cdot\na)y_e+(Y\cdot\na)Y=\na p}\vsp &\inr\
(0,\9)\times\calo,\vsp \na\cdot\ Y=0\hspace*{4cm}&\inr\
(0,\9)\times\calo,\vsp Y=u&\onr\ (0,\9)\times\pp\calo,\vsp
Y(0)=y_0-y_e&\inr\ \calo.\earr\end{equation} Our main concern here is
the design of an oblique boundary feedback controller which
stabilizes exponentially the equilibrium state $y_e$, or,
equivalently, the zero solution to \eqref{e1.1prim}. The main step
toward this end is the stabilization of the linear system
corresponding to \eqref{e1.1prim} or, more generally, of the
Oseen--Stokes system
\begin{equation}\label{e1.2}
\barr{ll} \derp Yt-\nu\D Y+(Y\cdot\na)a+(b\cdot\na)Y=\na p&\inr\
(0,\9)\times\calo,\vsp \na\cdot\ Y=0&\inr\ (0,\9)\times\calo,\vsp
Y=u&\onr\ (0,\9)\times\pp\calo,\earr\end{equation}where $a,b\in
(C^2(\ov\calo))^d,\ \na\cdot a=\na\cdot b=0$ in $\calo$.
 Besides its significance as first order linear approximation of \eqref{e1.1prim}, this system models the dynamics of a Stokes flow with inclusion of a convection acceleration $(b\cdot\na)Y$ and also the disturbance flow induced by a moving body in a Stokes fluid flow.

In its complex form, the main result  of this work, Theorem \ref{t1}, amounts to
saying that there is a boundary feedback controller of the form
\begin{equation}\label{e1.3}
\barr{ll} u(t,x)=\eta \dd\sum^N_{j=1}\mu_j\(\int_\calo
Y(t,x)\ov\vf^*_j(x)dx\)(\phi_j(x)+\a(x)\vec n(x)),\\\hfill t\ge0,\
x\in\pp\calo,\earr\end{equation} which stabilizes exponentially
system \eqref{e1.1}. Here $\phi_j={\rm lin}\left\{\derp{\vf^*_i}n\right\}^N_{i=1}$ and $\{\vf^*_j\}^N_{j=1}$ is an
eigenfunction system for the adjoint of the Stokes--Oseen operator
\begin{equation}\label{e1.4a}
\barr{rcll} \call\vf&=&-\nu\D\vf+(a\cdot\na)\vf+(\vf\cdot\na)b-\na
p,\ \vf\in D(\call),\vsp D(\call)&=&\{\vf\in(H^2(\calo))^d\cap
(H^1_0(\calo))^d;\ \na\cdot\vf=0\ \inr\
\calo\}.\earr\end{equation} It turns out (see Theorem \ref{t3})
that this feedback controller also stabilizes the \NS\ system
\eqref{e1.1prim} in a \nei\ of the origin.

In \eqref{e1.3}, $N$ is the number of the eigenvalue $\lbb_j$ of
$\call$ with $\Res\,\lbb_j<0$, $\a\in C(\ov\calo)$ is an
arbitrary function with zero circulation on $\pp\calo$, that is,
\begin{equation}\label{e1.4}
\barr{ll}
\dd\int_\calo\a(x) dx=0 \earr
\end{equation}
$\mbox{and }\phi_j={\rm lin\ span}\left\{\derp{\vf^*_j}n\right\}^N_{i=1}.$ (In its real version given in Theorem \ref{t4}, the stability controller $u$ is of the form \eqref{e1.3} but with $\{\vf^*_j\}$ replaced by $\{\Res\,\vf_j,\Ims\,\vf_j\}$.)

Taking into account that the vectors $\derp{\vf^*_j}n\ (x)$ are tangential at each $x\in\pp\calo$, that is, $\derp{\vf^*_j}\nu=\nu\vec\tau$, where $\vec\tau$ is the tangent vector (see \cite{11}, p.~35), we see that $u$ is an oblique vector field on $\pp\calo$.

More precisely, we have
\begin{eqnarray}
&u(t,x)\cdot\vec n(x)=\a(x),\ \ \ff x\in\pp\calo,\label{e1.5}\\[2mm]
&|\cos\<u(t,x),\vec n(x)\>|\ge 1-\dd\frac C{C+|\a(x)|}\,,\ \ \ff x\in\pp\calo,\label{e1.6}
\end{eqnarray}
where $C>0$ is independent of $\a$. This means that the "stabilizable boundary controller $u$ can be chosen "almost" normal to $\pp\calo$. However, for technical reasons the limit case $|\a|\equiv+\9$, that is, $u$ normal is excluded from our discussion.

It should be said that in the stabilization literature, only in a few situations was designed a normal stabilizable controller for \eq\ \eqref{e1.1} and this for periodic flows in $2{-}D$ channels (see, e.g., \cite{1}, \cite{2}, \cite{3}, \cite{20}, \cite{21}, \cite{22}). However, even in this case, the feedback controller is not given in explicit form and sometimes one assumes restrictive conditions on $\nu$  or on the spectrum of the operator $\call$.

It should be said that there is a large body of results obtained in recent years on boundary stabilization of system \eqref{e1.1} and here the works \cite{10}, \cite{11}, \cite{12}, \cite{13}, \cite{15}, \cite{16}, \cite{17}, \cite{18}  should be primarily cited. (See, also, \cite{7}, \cite{13}, \cite{14}, \cite{20a}.) The approach used in these works can be described in a few words as follows; one decomposes  system \eqref{e1.1} in a finite-dimensional unstable part which is exactly controllable and an infinite-dimensional part which is exponentially stable and proves so its stabilization by open loop boundary controller with finite-dimensional structure. Then one designs in a standard way a stabilizable feedback controller via the infinite-dimensional algebraic Riccati equation associated with an infinite horizon quadratic optimal control problem. Our construction of boundary stabilizable controller for \eqref{e1.1} avoids the Riccati \eq\ based approach which though provides a robust controller it is, however,  untreatable from computational point of view. Instead, we  propose an explicit feedback controller of the form \eqref{e1.3} easy to implement into system. It should be said that this construction resembles the form of stabilizable noise controllers recently designed in the author's works \cite{4}, \cite{5}, \cite{6}, \cite{8}, \cite{9}, which seem to be, however, more robust to stochastic perturbations.

The plan of the paper is the following. In Section 2, we present the main stabilization result which will be proved in Section 3. In Section 4, we shall give an application to stabilization of Stokes--Oseen periodic flows in a $2-D$ channel.

Everywhere in the following, we shall use the standard notation for spaces of functions on $\calo\3\rr^d$. In particular, $C^k(\ov\calo)$, $k=0,1,...$, is the space of $k$-differentiable functions on $\ov\calo$ and $H^k(\calo)$, $k=1,2,$ $H^1_0(\calo)$ are Sobolev spaces on $\calo$.

\section{The main result}
\setcounter{equation}{0}

\subsection{Notation}

Everywhere in the following, $\calo$ is a bounded and open domain of $\rr^d$, $d=2,3$, with smooth and simply connected boundary $\pp\calo$.

We set
$$H=\{y\in(L^2(\calo))^d;\ \na\cdot y=0\ \inr\ \calo,\ \ y\cdot\vec n=0\ \onr\ \pp\calo\}$$and denote by $\Pi:(L^2(\calo))^d\to H$ the Leray projector on $H$. We consider the operator $A:D(A)\3H\to H$, $\cala:D(\cala)\3H\to H$,
\begin{eqnarray}
\quad Ay&=&-\nu\Pi(\D y),\ \ff y\in D(A)=(H^1_0(\calo))^d\cap(H^2(\calo))^d\cap H,\label{e2.1}\\[2mm]
\quad Ay&=&\Pi(-\nu\D y+(y\cdot\na)a+(b\cdot\na)y)\label{e2.2}\\[1mm]
&=&Ay+\Pi((y\cdot\na)a+(b\cdot\na)y),\nonumber\\[1mm] &&\qquad\qquad\ff y\in D(\cala)=D(A).\nonumber
\end{eqnarray}
We denote by $\wt H$ the complexified space $\wt H=H+iH$ and consider the extension $\wt\cala$ of $\cala$ to $\wt H$, that is, $\wt\cala(y+iz)=\cala y+i\cala z$ for all $y,z\in D(\cala)$.

The scalar product of $H$ and of $\wt H$ are denoted by $\<\cdot,\cdot\>$ and $\<\cdot,\cdot\>_{\wt H}$, respectively. The corresponding norms are denoted by $|\cdot|_H$ and $|\cdot|_{\wt H}$, respectively.

For simplicity, we denote in the following again by $\cala$ the operator $\wt\cala$ and the difference will be clear from the content. The operator $\cala$ has a compact resolvent $(\lbb I-\cala)^{-1}$ (see, e.g., \cite{7}, p~92). Consequently, $\cala$ has a countable number of eigenvalues $\{\lbb_j\}^\9_{j=1}$   with corresponding eigenfunctions $\vf_j$ each with finite algebraic multiplicity $m_j$. In the following, each eigenvalue $\lbb_j$ is  repeated according to its algebraic multiplicity $m_j$.

Note also that   there  is a finite number of eigenvalues $\{\lbb_j\}^N_{j=1}$ with $\Res\,\lbb_j{\le} 0$ and that the spaces $X_u={\rm lin\ span}\{\vf_j\}^N_{j=1}=P_N\wt H,$ $X_s=(I-P_N)\wt H$  are invariant with respect to $\cala$. Here, $P_N$  is the algebraic projection of $\wt H$ on $X_u$ and is defined by
\begin{equation}\label{e2.2a}
P_N=\frac1{2\pi i}\int_\G(\lbb I-\cala)\1 d\lbb,\end{equation}where $\G$ is a closed curve which contains in  interior the eigenvalues $\{\lbb_j\}^N_{j=1}$.

If we set $\cala_u=\cala|_{X_u},\ \cala_s=\cala|_{X_s}$, then we have
$$\sigma(\cala_u)=\{\lbb_j:\Res\,\lbb_j\le0\},\ \ \sigma(\cala_s)=\{\lbb_j:\Res\,\lbb_j>0\}.$$

We recall that the eigenvalue $\lbb_j$ is called semisimple if its algebraic multiplicity  $m_j$ coincides with its geometric multiplicity $m^g_j$.   In particular, this happens if $\lbb_j$ is simple and it turns out that the property of the eigenvalues $\lbb_j$ to be all simple is generic (see \cite{7}, p.~164). The dual operator $\cala^*$ has the eigenvalues $\ov\lbb_j$ with the eigenfunctions $\vf^*_j$, $j=1,...\,.$

For the time being, the following hypotheses will be assumed.
\bit\item[(H1)] {\it The eigenvalues $\lbb_j,\ j=1,...,N,$ are semisimple.}\eit
This implies that
\begin{equation}\label{e2.3}
\cala\vf_j=\lbb_j\vf_j,\ \ \cala^*\vf^*_j=\ov\lbb_j\vf^*_j,\ \ \ j=1,...,N,\end{equation}
and so we can choose systems $\{\vf_j\},\{\vf^*_j\}$ in such a way that
\begin{equation}\label{e2.4}
\<\vf_j,\vf^*_k\>_{\wt H}=\de_{jk},\ \ j,k=1,...,N.\end{equation}
Next hypothesis is a unique continuation assumption on normal derivatives $\derp{\vf^*_j}n\,,$ $j=1,...,N.$
\bit\item[(H2)] {\it The system $\left\{\derp{\vf^*_j}n\right\}^N_{j=1}$ is linearly independent on $\pp\calo$.}\eit
We note that,  in the special case $N=1$, hypothesis (H2) reduces to: $\derp{\vf^*_j}n\not\equiv0$ for all $j=1,...,N.$  It is not known if this unique continuation property is always satisfied,    but it holds, however, for "almost all $a,b$" in the generic sense (see \cite{12a}). In specific examples, however, this assumption might be easily checked and we shall see later on in Section 4 that it holds for systems in a $2-D$ channel $\calo=\{(x,y)\in\rr\times(0,1)\}$ with periodic conditions in $x$.

\subsection{The main stabilization result}

Consider the feedback boundary controller
\begin{equation}\label{e2.5}
u=\eta \sum^N_{j=1}\mu_j\<P_NY,\vf^*_j\>_{\wt H}(\phi_j+\a\vec n),\end{equation}where
\begin{eqnarray}
\mu_j&=&\frac{k+\lbb_j}{k+\lbb_j-\nu\eta}\,\ \ j=1,...,N,\label{e2.6}\\[2mm]
\phi_j&=&\sum^N_{i=1}\a_{ij}\ \derp{\vf^*_i}n\,,\ \ j=1,...,N,\label{e2.7}
\end{eqnarray}
and the matrix $\calx=\|\a_{ij}\|^N_{i,j=1}$ is given by
\begin{equation}\label{e2.8a}
\calx=\calf\1,\ \calf=\left\|\int_{\pp\calo}\derp{\vf^*_i}n\cdot\derp{\ov\vf^*_j}n\ dx\right\|^N_{i,j=1}.\end{equation}
In virtue of hypothesis (H2), $\calf$ is invertible and so $\calx$ is well defined.

\begin{theorem}\label{t1} Assume that {\rm(H1), (H2),
\eqref{e1.4}} hold and that $\Res\,\lbb_j<0$ for $j=1,...,N$,
$\Res\ \lbb_j>0$ for $j>N.$ Let $k>0$ sufficiently large and
$\eta>0$ be such that
\begin{equation}\label{e2.8}
(|k+\lbb_j|^2-\eta k\nu)\Res\, \lbb_j-\eta\nu\ \Res\,\lbb_j^2>0.\end{equation}Then the feedback controller \eqref{e2.5}
stabilizes exponentially system \eqref{e1.2}, that is, the
solution $Y$ to the closed loop system
\begin{equation}\label{e2.9}
\barr{ll} \derp Yt-\nu\D Y+(Y\cdot\na)a+(b\cdot\na)Y=\na p&\inr\
(0,\9)\times\calo,\vsp \na\cdot Y=0&\inr\ (0,\9)\times\calo,\vsp
Y=\eta\dd\sum^N_{j=1}\mu_j\<P_NY,\vf^*_j\>_{\wt H}(\phi_j+\a\vec
n)&\onr\ (0,\9)\times\pp\calo,\earr\end{equation} satisfies
\begin{equation}\label{e2.10}
|Y(t)|_{\wt H}\le Ce^{-\g t}|Y(0)|_{\wt H},\ \ \ff
t\ge0,\end{equation}for some $\g>0.$\end{theorem}

As noticed earlier, by \eqref{e2.7} it follows that, in each $x\in \pp\calo$,  $\phi_j(x)$ are tangent to $\pp\calo$ and so, for $|\a|$ large enough, the controller $u$ is "almost" normal. Moreover, since $\Res\,\lbb_j<0$ for $j=1,...,N$, by \eqref{e2.6} it is easily seen that \eqref{e2.8} holds for $\eta>0$  and $k>0$ sufficiently large and suitable chosen.

It should be observed that, if assumption (H2) is strengthen to all $j=1,...,$ and so \eqref{e2.4} holds for all $i,j=1,...,$ then $\<P_NY,\vf^*_j\>_{\wt H}=\<\psi,\vf^*_j\>_{\wt H}$ for all $j$ and so the controller \eqref{e2.5} reduces to
$$u=\eta\sum^N_{j=1}\mu_j\<Y,\vf^*_j\>_{\wt H}(\phi_j+\a\vec n).$$
If $\lbb_j$ are complex valued, then the controller \eqref{e2.5} is complex valued too and plugged into system \eqref{e1.2} leads to a real closed loop system in $(\Res\,Y,\Ims\,Y)$. In order to avoid this situation, we shall construct in Section 3.3 a real stabilizable feedback controller of the form \eqref{e2.5} which has a similar stabilization effect. (See Theorem \ref{t4}.)

\begin{remark}\label{r2.1} {\rm With choice of $\phi_j$ we have $\dd\int_{\pp\calo}\ov\phi_j\ \derp{\vf^*_j}n\ dx=\de_{ij},\ i,j=1,...,N,$ and, as seen  later on, this is essential in the proof of Theorem \ref{t1}. However, this can be also achieved for $\phi_j$ of the form
$$\phi_j=\dd\sum^N_{i=1}\alpha_{ij}\chi_i$$
where $\{\chi_i\}$ are suitably chosen.

To find such $\chi_i$ and $\alpha_{ij}$, it sufficed to   assume instead (H2) that all $\derp{\vf^*_j}n\not\equiv0$ on $\pp\calo$.}\end{remark}

\subsection{Stabilizable controllers with support in $\G_0\subset\pp\calo$}

Consider system \eqref{e1.1} with a boundary controller $u$ with support in an open and smooth subset $\G_0\3\pp\calo$, that is,
\begin{equation}\label{e2.11}
\barr{ll}
\derp Yt-\nu\D Y+(Y\cdot\na)a+(a\cdot\na)Y=\na p&\inr\ (0,\9)\times\calo,\vsp
\na\cdot Y=0\ \ \inr\ (0,\9)\times\calo,\vsp
Y=\one_{\G_0}u&\onr\ (0,\9)\times\pp\calo,\earr\end{equation}
where $\one_{\G_0}$ is the characteristic function of $\G_0$.

In this case, instead of (H2) we assume that
\bit\item[(H2)$'$] {\it The system $\left\{\derp{\vf^*_j}n\right\}^N_{j=1}$ is linearly independent on $\G_0$.}\eit
We assume also  that
\begin{equation}\label{e2.12}
\int_{\G_0}\a(x)dx=0.\end{equation}We choose $\wt\phi_j,\ j=1,...,N$, of the form
\begin{equation}\label{e2.13}
\wt\phi_j=\sum^N_{k=1}\wt\a_{jk}\ \derp{\vf^*_j}n\,,\end{equation}where the matrix $\|\wt\alpha_{jk}\|^N_{j,k=1}$ is given by
$$\left(\left\|\int_{\pp\calo}\derpp n\ \vf^*_j\cdot\derpp n\ \vf^*_k\, dx\right\|^N_{i,j=1}\right)\1.$$

Consider the feedback controller
\begin{equation}\label{e2.14}
u_{\G_0}=\eta\sum^N_{j=1}\mu_j\<P_NY,\vf^*_j\>_{\wt H}(\wt\phi_j+\a\vec n),\end{equation}where $\eta,\mu_j$ are chosen as in Theorem \ref{t1}.

We have

\begin{theorem}\label{t2} The controller $u_{\G_0}$ stabilizes
exponentially system \eqref{e1.2}.\end{theorem}

\subsection{Stabilization of system \eqref{e1.1}
(\eqref{e1.1prim})}

We set $W=(H^{\frac12-\vp}(\calo))^d\cap H$ for
$Z=(H^{\frac32+\vp}(\calo))^d\cap H)$, $d=2.$

\begin{theorem}\label{t3} Let $d=2$. Then, under the assumptions
of Theorem {\rm\ref{t1}}, the feedback boundary controller
\eqref{e2.5} stabilizes exponentially system \eqref{e1.1prim} in a
\nei\ $\calw=\{y_0\in W;\ \|y_0\|_W<\rho\}$. More precisely, the
solution $y\in C([0,\9);W)\cap L^2(0,\9;Z)$ to the closed loop
system
\begin{equation}\label{e2.15}
\barr{ll}
\derp Yt-\nu\D Y+(Y\cdot\na)a+(b\cdot\na)Y+(Y\cdot\na)Y=\na p\\\hfill \inr\ (0,\9)\times\calo,\vsp
Y=\eta\dd\sum^N_{j=1}\mu_j\<P_NY,\vf^*_j\>_{\wt H}(\phi_j+\a\vec n)\ \onr\ (0,\9)\times\pp\calo,
\earr\end{equation}satisfies for $Y(0)\in \calw$ and
$\rho$ sufficiently small
\begin{equation}\label{e2.16}
\|Y(t)\|_W\le C e^{-\g t}\|Y(0)\|_W,\ \ \ff t\ge0,\end{equation}for some
$\g>0$.
\end{theorem}
In particular, it follows that the boundary feedback controller
\begin{equation}\label{e2.17}
u=\eta\sum^N_{j=1}\mu_j\<P_N(Y-y_e),\vf^*_j\>_{\wt
H}(\phi_j+\a\vec n)\end{equation}stabilizes exponentially the
equilibrium solution $y_e$ to \eqref{e1.1} in a \nei\ $\{y_0\in
W;\ \|y_0-y_e\|_W<\rho\}$.

\section{Proofs}
\setcounter{equation}{0}

\subsection{Proof of Theorem \ref{t1}}

We set
$$U^0=\lac u\in(L^2(\pp\calo))^d;\int_{\pp\calo}u(x)\cdot\vec n(x)dx=0\rac.$$
Then, for $k>0$ sufficiently large, there is a unique solution $y\in (H^{\frac12}(\calo))^d$ to the \eq
$$\barr{l}
-\nu\D y+(y\cdot\na)a+(b\cdot\na)y+ky=\na p\ \inr\ \calo,\vsp
\na\cdot y=0\ \inr\ \calo,\ \ y=u\ \onr\ \pp\calo.\earr$$
(See, e.g., \cite{19a}, p. 365.) We set $y=Du$ and note that (see, e.g., \cite{11}, p.~102),
$$D\in L((H^s(\pp\calo))^d\cap U^0;(H^{s+\frac12})\calo))^d),\mbox{ for }s\ge-\frac12\,\cdot$$In terms of the Dirichlet map $D$, system\eqref{e1.2} can be written as
\begin{equation}\label{e3.1}
\barr{l}
\dd\frac d{dt}\ Y(t)+\cala(Y(t)-Du(t))=0,\ \ t\ge0,\vsp
Y(0)=y_0.\earr\end{equation}Equivalently,
\begin{equation}\label{e3.2}
\barr{l}
\dd\frac d{dt}\ z(t)+\cala z(t)=-\Pi\(D\ \dd\frac{du}{dt}\ (t)\),\ \ t\ge0,\vsp
z(0)=y_0-Du(0),\earr
\end{equation}
\begin{equation}\label{e3.3}
z(t)=Y(t)-Du(t),\ \ t\ge0.\end{equation}
In the \fwg, we fix $k>0$ sufficiently large and $\eta>0$ such that \eqref{e2.8} holds. In particular, we also have
\begin{equation}\label{e3.4}
 \lbb_i+k-\nu\eta\ne0\ \forr\ i=1,2,...,N.\end{equation}

We note fist that in terms of $z$ the controller \eqref{e2.5} can be, equivalently, expressed as
\begin{equation}\label{e3.5}
u(t)=\eta\sum^N_{j=1}\<P_N z(t),\vf^*_j\>_{\wt H}(\phi_j+\a\vec n).\end{equation}Indeed, by \eqref{e3.3} and \eqref{e3.5}, we have
\begin{equation}\label{e3.6}
\barr{lcl}
u(t)&=&\dd\eta\sum^N_{j=1}\<P_NY(t),\vf^*_j\>)_{\wt H}(\phi_j+\a\vec n)\vsp&&-\eta\dd\sum^N_{j=1}\<u(t),D^*\vf^*_j\>_{(L^2(\pp\calo))^d}(\phi_j+\a\vec n),\earr\end{equation}where $D^*$ is the adjoint of $D$.

On the other hand, if we set $\psi=D(\phi_j+\a\vec n)$ and recall that
$$\barr{rclcrcl}
\call^*\vf^*_i-\ov\lbb_i\vf^*_i&=&\na p_i\ \inr\ \calo,&&\vf^*_i&=&0\ \onr\ \pp\calo,\vsp
\call\psi+k\psi&=&\na\wt p\ \inr\ \calo,&&\psi&=&\phi_j+\a\vec n\ \onr\ \pp\calo,\earr$$where $\call$ is the Stokes--Oseen operator \eqref{e1.4a} and $\call^*$ is its formal adjoint, we get via Green's formula
\begin{equation}\label{e3.7}
\barr{r}\<\phi_j{+}\a\vec n,D^*\vf^*_i\>_{(L^2(\pp\calo))^d}{=}\!\dd\int_\calo\!\!
\psi\cdot\ov\vf^*_i dx=-\dd\frac\nu{\lbb_i{+}k}\!\int_{\pp\calo}\!\!(\phi_j{+}\a\vec n){\cdot}\derp{\ov\vf^*_i}n\ dx\vsp\qquad\dd=-\dd\frac\nu{\lbb_i{+}k}\ \de_{ij},\ \ff i,j=1,...,N,\earr\end{equation}because $\vec n\cdot\derp{\vf^*_i}n=0,$ a.e. on $\pp\calo$ (see \cite{11}, Lemma 3.3) and, by \eqref{e2.4}, \eqref{e2.7}, \eqref{e2.8a}, we have
\begin{equation}\label{e3.8}
\int_{\pp\calo}\phi_j\cdot\derp{\ov\vf^*_i}n\ dx=\de_{ij},\ \ i,j=1,...,N.\end{equation}Then, by \eqref{e3.6}, \eqref{e3.7}, we see that
$$\<u(t),D^*\vf^*_i\>_{(L^2(\pp\calo))^d}=\dd\frac{-\eta\nu}
{k+\lbb_i-\nu\eta}\ \<P_NY,\vf^*_i\>_{\wt H}$$and, substituting into \eqref{e3.6}, we get \eqref{e2.5} as claimed.

Now, substituting \eqref{e3.5} into \eqref{e3.2}, we obtain that
\begin{equation}\label{e3.9}
\barr{l}
\dd\frac{dz}{dt}+\cala z=-\eta\dd\sum^N_{j=1}\<P_N\ \dd\frac d{dt} z(t),\vf^*_j\>\Pi D(\phi_j+\a\vec n),\vsp
z(0)=z_0=y_0-Du(0).\earr\end{equation}

We write \eqref{e3.9} as
\begin{eqnarray}
\quad\frac{dz_u}{dt}+\cala_u z_u&=&-\eta P_N\dd\sum^N_{j=1}\<P_N\dd\frac{dz}{dt}\,,\vf^*_j\>_{\wt H}\Pi D(\phi_j+\a\vec n),\label{e3.10}\\[2mm]
\quad\frac{dz_s}{dt}+\cala_s z_s&=&-\eta (I-P_N)\dd\sum^N_{j=1}\<P_N\dd\frac{dz}{dt}\,,\vf^*_j\>_{\wt H}\Pi D(\phi_j+\a\vec n),\label{e3.11}
\end{eqnarray}
$z=z_u+z_s,\ z_u\in X_u,$ $z_s\in X_s$ and $P_N$ is given by \eqref{e2.2a}. If we represent $z_u$ as $$z_u=\dd\sum^N_{j=1}z_j\vf_j,$$
and recall \eqref{e3.7}, we rewrite \eqref{e3.10} as
\begin{equation}\label{e3.12}
z'_j+\lbb_jz_j=\frac{\eta\nu}{k+\lbb_j}\ z'_j,\ \ t\ge0.\end{equation}
By \eqref{e2.8} we have
$$\Res\left[\lbb_j
\(1-\frac{\eta\nu}{k+\lbb_j}\)^{-1}\right]>0.$$
Then, by \eqref{e3.12} we see that  we have for some $\g_0>0$
\begin{equation}\label{e3.13}
|z_j(t)|\le e^{-\g_0t}|z_j(0)|,\ \ j=1,...,N.\end{equation}On the other hand, by \eqref{e3.11} we have
\begin{equation}\label{e3.14}
\frac{dz_s}{dt}+\cala_s z_s=\eta(I-P_N)\sum^N_{j=1}z_j\Pi D(\phi_j+\a\vec n),\end{equation}and since
$$\|e^{-\cala_st}\|_{L(\wt H,\wt H)}\le Ce^{-\g_1t},\ \ \ff t\ge0,$$
for some  $\g_1>0$, we see that
$$|z_s(t)|_{\wt H}\le C\exp(-\g_0t)|z_s(0)|_{\wt H},\ \ \ff t\ge0,$$which together with \eqref{e3.13} yields
\begin{equation}\label{e3.15}
|z(t)|_{\wt H}\le C\exp(-\g_0 t)|z(0)|_{\wt H},\ \ \ff t\ge0.\end{equation}

Now, recalling \eqref{e3.3} and \eqref{e3.5}, we obtain \eqref{e2.10}, thereby completing the proof.

\subsection{Proof of Theorem \ref{t2}}

The proof is exactly the same as that of Theorem \ref{t1} except that the Dirichlet map $D$ is taken for the boundary condition $y=\one_{\G_0}.$ The details are   omitted.

\subsection{Proof of Theorem \ref{t3}}

We shall apply Theorem 1.2.1 from \cite{10} (see, also, Theorem 5.1 in \cite{11}).

In fact, system \eqref{e1.1prim} with the feedback controller $$u=FY=\eta\dd\sum^N_{j=1}\mu_j\<P_NY,\vf^*_j\>_{\wt H}(\phi_j+\a\vec n)$$ can be written as
$$\barr{l}
\dd\frac{dY}{dt}+\cala(Y-DFY)+BY=0,\ t>0,\vsp
Y(0)=y_0,\earr$$where $BY=\Pi(Y\cdot\na)Y).$

By Theorem \ref{t1}, it is easily seen that the operator $A_F=\cala(I-DF):W\to W$ with $D(A_F)=\{y\in W;\ \cala(y-DFy)\in W\}$ generates an analytic $C_0$-semigroup on $W$ which is exponentially stable on $W$.

Moreover, coming back to system \eqref{e3.9}-\eqref{e3.11}, we see that besides \eqref{e3.15} we have also
$$\int^\9_0|A^{\frac34}z(t|^2dt\le C\|z(0)\|^2_W$$and recalling that $Y=e^{-A_Ft}y_0$ is given by
$$Y=z+DFY=z+\eta\sum^N_{j=1}\mu_j\<P_NY,\vf^*_j\>_{\wt H}D(\phi_j+\a\vec n)$$we infer that
$$\int^\9_0\|e^{-A_Ft}y_0\|^2_Zdt\le C\|y_0\|^2_W,\ \ \ff y_0\in W.$$Then, by Theorem 1.2.1 from \cite{10}, we infer that the conclusion of Theorem \ref{t3} holds.

\subsection{Real stabilizable feedback controllers}

We shall construct here a real stabilizable feedback controller of the form \eqref{e3.5}.
To this purpose, we consider in the space $H$ the system $$\{\Res\,\vf_j,\Ims\,\vf_j\}^N_{j=1}=\{\psi_j\}^N_{j=1}.$$ We set $X^*_u={\rm lin\ span}\{\Res\,\vf_j,\Ims\ \vf_j\}^N_{j=1}$, $j=1,...,N.$ We decompose the space $H=X^*_u\oplus X^*_s$ and note that the real operator $\cala$ leaves invariant both spaces $X^*_s$ and $X^*_u$ and $\cala^*_s=\cala|_{X^*_s}$ generates an exponential stable semigroup on $X^*_s\3H$.

We have
\begin{equation}\label{e3.16a}
\cala\psi_j=(\Res\,\lbb_j)\psi_j-(\Ims\,\lbb_j)\psi_{j+1},\
\cala\psi_{j+1}=(\Ims\,\lbb_j)
\psi_j+(\Res\,\lbb_j)\psi_{j+1}.\end{equation}
We may assume via Schmidt's ortogonalization algorithm that the system $\{\psi_j\}^N_{j=1}$ is orthonormal.

Then, we construct the feedback controller
\begin{equation}\label{e3.16}
u^*=\eta\sum^N_{j=1}\mu_j\<P_NY,\psi_j\>(\phi^*_j+\a\vec n),\end{equation}where $\phi^*_j$ is of the form
\begin{equation}\label{e3.17}
\phi^*_j=\sum^N_{i=1}\a^*_{ij}\ \derp{\psi_i}n\,,\ \ j=1,...,N,\end{equation}and $\a^*_{ij}$ are chosen in a such a way that
$$\sum^N_{i=1}\a^*_{ij}\<\wt D\ \derp{\psi_i}n\,,\psi_\ell\> =\de_{j\ell},\ j,\ell=1,...,N,$$
where $\wt D$ is the Dirichlet map corresponding to the operator $\cala^*_k$. Keeping in mind  that
$$\<\wt D\chi_i,\cala_k\psi_j\> =-\nu\int_{\pp\calo}\chi_i\ \derpp n\ \psi_j\, dx,\ i,j=1,...,N,$$
we see  by \eqref{e3.16a}  that, for $k$ large enough, we have for $\chi_i=\derp{\psi_i}n\,,$
$$\barr{rcl}
\<\wt D\chi_i,\psi_j\>&=&-\dd\frac{\nu(b^i_j\ \Res\,\lbb_j+b^i_{j+1}\ \Ims\,\lbb_j)}{|z_j|^2+k\ \Res\,\lbb_j},\vsp
\<\wt D\chi_i,\psi_{j+1}\>&=&-\dd\frac{\nu((\Res\,\lbb_j+k)b^i_{j+1}-\Ims\,\lbb_j b^i_j)}{|z_j|^2+k\ \Res\,\lbb_j},\earr$$
where $b^i_j=\dd\int_{\pp\calo}\chi_i\ \derp{\psi_j}n\ dx$.
Then, assuming that
\bit\item[(H2)$^*$] {\it The system $\left\{\derp{\vf_j}n\right\}^N_{j=1}$ is linearly independent on $\pp\ooo$.}\eit
it follows that so is $\left\{\derp{\psi_j}n\right\}^N_{j=1}$ and this implies that such a choice of $\a^*_{ij}$   is possible. Then, arguing exactly as in the proof of Theorem \ref{t1}, we see that, for $\eta$ and $\mu_j$ suitable chosen, the real controller \eqref{e3.16} stabilizes exponentially system \ref{e1.2}. We have, therefore,

\begin{theorem}\label{t4} Under  assumptions {\rm(H1), (H2)$^*$} and
\eqref{e1.4}, there is a boun\-dary feedback controller $u^*$ of
the form \eqref{e3.16} which stabilizes exponentially
system~\eqref{e1.2}.\end{theorem}

The proof is exactly the same as that of Theorem \ref{t1} and so it is omitted.

We note, however, that if instead of (H2)$^*$ we assume only that $\derp{\vf_j}n\not\equiv0$ on $\pp\calo$ for $j=1,...,N,$ then Theorem \ref{t4} still remains valid with $\phi^*=\dd\sum^N_{i=1}\alpha_{ij}\chi_i$, where $\chi_i$ are chosen in such a way that $\dd\sum^N_{i=1}\alpha^*_{ij}\<\wt D\chi_i,\psi_\ell\>=\delta_{j\ell}$, $j,\ell=1,...,N.$ Note also that Theorem \ref{t3} remains true in the present situation.

\section{Boundary stabilization of a periodic flow\\ in a $2{-}D$
channel} \setcounter{equation}{0}

Consider a laminar flow in a two-dimensional channel with the walls located at $y=0,1.$ We shall assume that the velocity field $(u(t,x,y),v(t,x,y))$ and the pressure $p(t,x,y)$ are $2\pi$ periodic in $x\in(-\9,+\9)$.

The dynamic of flow is governed by the incompressible $2-D$ \NS\ \eqs
\begin{equation}\label{e4.1}
\barr{ll}
u_t-\nu\D u+uu_x+vu_y=p_x,\quad x\in\rr,\ y\in(0,1),\vsp
v_t-\nu\D v+uv_x+vv_y=p_y,\quad x\in\rr,\ y\in(0,1),\vsp
u_x+v_y=0,\vsp
u(t,x+2\pi,y)\equiv u(t,x,y),\ \ v(t,x+2\pi,y)\equiv v(t,x,y),\ \ y\in(0,1).\earr\end{equation}

Consider a steady-state flow governed by \eqref{e4.1} with zero vertical ve\-lo\-city com\-po\-nent, i.e., $(U(x,y),0)$. Since the flow is freely divergent, we have \mbox{$U_x\equiv0
$} and so $U(x,y)\equiv U(y)$. This yields
\begin{equation}\label{e4.2}
U(y)=C(y^2-y),\ \ \ff y\in(0,1),\end{equation}where $C\in\rr^-$. In the \fwg, we take $C=-\dd\frac a{2\nu}$ where $a\in\rr^+$.

The linearization of \eqref{e4.1} around the steady-state flow $(U(y),0)$ leads to the \fwg\ system
\begin{equation}\label{e4.3}
\barr{ll}
u_t-\nu\D u+u_xU+vU'=p_x,\ \ y\in(0,1),\ x,t\in\rr,\vsp
v_t-\nu\D v+v_xU=p_y,\vsp
u_x+v_y=0,\vsp
u(t,x+2\pi,y)\equiv u(t,x,y),\ \ v(t,x+2\pi,y)\equiv v(t,x,y).\earr\end{equation}
Here  we apply Theorem \ref{t1} to construct an oblique boundary feedback controller for system\eqref{e4.3}. To this aim, we recall first the Fourier functional setting for description of periodic fluid flows in the channel $(-\9,+\9)\times(0,1).$

Let $L^2_\pi(Q)$, $Q=(0,2\pi)\times(0,1)$ be the space of all the functions $u\in L^2_{\rm loc}(R\times(0,1))$ which are $2\pi$-periodic in $x$. These functions are characterized by their Fourier series

 $$\mbox{$u(x,y)=\dd\sum_ka_k(y)e^{ikx},\ a_k=\bar a_{-k},$ $a_0=0$, $\dd\sum_k\int^1_0|a_k|^2dy<\9.$}$$

Similarly, $H^1_\pi(Q),\ H^2_\pi(Q)$ are defined. For instance,
$$\barr{l}
H^1_\pi(Q)=\lac u\in L^2_\pi(Q);\ u\in\dd\sum_ka_ke^{ikx},\ a_k=\bar a_{-k},\ a_0=0,\right.\vsp \qquad\qquad\qquad\left.\dd\sum_k \int^1_0(k^2|a_k|^2+|a'_k|^2)dy<\9\rac,\ k\to j.\earr$$We set
$$H=\{(u,v)\in(L^2_\pi(Q))^2;\ u_x+v_y=0,\ v(x,0)=v(x,1)=0\}.$$If $u_x+v_y=0$, then the trace of $(u,v)$ at $y=0,1$ is well defined as an element of $H\1(0,2\pi)\times H\1(0,2\pi)$ (see, e.g., \cite{19}).

We also set
$$V=\{(u,v)\in H\cap H^1_\pi(Q);\ u(x,0)=u(x,1)=v(x,0)=v(x,1)=0\}.$$As defined above, the space $L^2_\pi(Q)$ is, in fact, the factor space $L^2_\pi(Q)/Z.$ The space $H$ can be defined equally as
$$\barr{lcl}H&=&\Big\{ u=\dd\sum_{k\ne0}u_k(y)e^{ikx},\ v=\dd\sum_{k\ne0}v_k(y)e^{ikx},\ v_kj(0)=v_k(1)=0, \vsp
&&\ \ \ \ \ \dd\sum_{k\ne0}\int^1_0(|u_k|^2+|v_k|^2)dy<\9,\ iku_k(y)+v'_k(y)=0,\vsp
&&\hfill\ \ \ \  \mbox{ a.e. }y\in(0,1),\ k\in\rr\Big\},\  k\to j.\earr$$

Let $\Pi:L^2_\pi(Q)\to H$ be the Leray projector and $\cala:D(\cala)\3H\to H'$ the operator
\begin{equation}\label{e4.4}
\barr{r}
\cala(u,v)=\Pi\{-\nu\D u+u_xU+vU',\ -\nu\D v+v_xU\},\vsp
\ff(u,v)\in D(\cala)=(H^2((0,2\pi)\times(0,1)).\earr
\end{equation}We associate with \eqref{e4.3} the boundary value conditions
\begin{equation}\label{e4.5}
\barr{lll}
u(t,x,0)=u^0(t,x),&u(t,x,1)=u^1(t,x),&t\ge0,\ x\in\rr,\vsp
v(t,x,0)=v^0(t,x),&v(t,x,1)=v^1(t,x),&t\ge0,\ x\in\rr,\earr
\end{equation}and, for $k>0$ sufficiently large, we consider the Dirichlet map $D:X\to L^2_\pi(Q)$ defined by $D(u^*,v^*)=(\wt u,\wt v)$,
\begin{equation}\label{e4.6}
\barr{ll}
-\nu\D\wt u+\wt u_\lbb U+\wt vU'+k\wt u=p_x,\ x\in\rr,\ y\in(0,1),\vsp
-\nu\D\wt v+\wt v_xU+k\wt v=p_y,\ x\in\rr,\ y\in(0,1),\vsp
\wt u_x+\wt v_y=0,\ \wt u(x+2\pi,y)=\wt u(x,y),\ \wt v(x+2\pi,y)=\wt v(x,y),\vsp
\wt u(x,y)=u^*(x,y),\ \wt v(x,y)=v^*(x,y),\ y=0,1.\earr\end{equation}Here
$$\barr{lcl}
X&=&\Big\{(u^*,v^*)\in L^2((0,2\pi)\times\pp(0,1));\ u^*(x+2\pi,y)=u^*(x,y), \vsp
&&\ \  v^*(x+2\pi,y)=v^*(x,y),\ \dd\int^{2\pi}_0 v^*(x,0)dx=\dd\int^{2\pi}_0 v^*(x,1)dx\Big\}.\earr$$Then system \eqref{e4.3} with boundary conditions \eqref{e4.4} can be written as
\begin{equation}\label{e4.7}
\barr{l}
\dd\frac d{dt}\ Y(t)+\cala(Y(t)-DU^*(t))=0,\ \ t\ge0,\vsp
Y(0)=(u_0,v_0),\earr\end{equation}where $Y=(u,v),\ U^*=(u^*,v^*).$

In order to apply Theorem \ref{t1}, we shall check hypothesis (H2) in this case.

To this end, we denote again by $\cala$ the extension of $\cala$ on the complexified space $\wt H$ and by $\lbb_j,\vf_j$ the eigenvalues and corresponding eigenvectors of the operator $\cala$. By $\vf^*_j$, we denote the eigenvector to the dual operator $\cala^*$.

\begin{lemma}\label{l1} For all $j=1,2,...,N,$ we have
\begin{equation}
\derp{\vf_j}n\ (x,y)\not\equiv0,\ \ x\in(0,2\pi),\
y=0,1,\label{e4.8}\end{equation}
and
\begin{equation}
\derp{\vf^*_j}n\ (x,y)\not\equiv0,\ \ x\in(0,2\pi),\
y=0,1.\label{e4.9}
\end{equation}\end{lemma}

\pf If we represent $\vf_j=(u^j,v^j),$ then \eqref{e4.8} reduces to
\begin{equation}\label{e4.10}
\left|\derpp y\ v^j(x,y)\right|+\left|\derpp y\ u^j(x,y)\right|>0,\ x\in(0,2\pi),\ y=0,1.\end{equation}
We set $\lbb=\lbb_j$ and $\vf_j=(u,v)$. This means that, if $\lbb$ is semisimple, then

\begin{equation}\label{e4.10a}
\barr{ll}
-\nu\D u+u_xU+vU'=\lbb u+p_x,\ x\in\rr,\ y\in(0,1),\vsp
-\nu\D v+v_xU=\lbb v+p_y,\ x\in\rr,\ y\in(0,1),\vsp
u_x+v_y=0,\vsp
u(x+2\pi,y)=u(x,y),\ v(x+2\pi,y)=v(x,y).\earr\end{equation}
If we represent $u,v,p$ as Fourier series,
\begin{equation}\label{e4.11}
u(x,y){=}\sum_k e^{ikx}u_k(y),\, v(x,y){=}\dd\sum_k e^{ikx}v_k(y),\, p(x,y){=}\dd\sum_ke^{ikx}p_k(y),\end{equation} we reduce \eqref{e4.10a} to (see, e.g., \cite{7}, p. 144)
$$\barr{l}
-\nu u''_k+(\nu k^2+ikU)u_k+U' v_k=ikp_k+\lbb u_k,\ y\in(0,1),\vsp
-\nu v''_k+(\nu k^2+ikU)v_k=p'_k+\lbb v_k,\ ik u_k+v'_k=0\ \inr\ (0,1),\vsp
u_k(0)=u_k(1)=0,\ \ v_k(0)=v_k(1)=0.\earr$$
Equivalently,
\begin{equation}\label{e4.12}
\barr{l}
-\nu v^{\rm iv}_k{+}(2\nu k^2{+}ikU)v''_k{-}k(\nu k^3{+}ik^2U{+}iU'')v_k{-}\lbb(v''_k{-}k^2v_k)=0,\\\hfill y\in(0,1),\vsp
v_k(0)=v_k(1)=0,\ v'_k(0)=v'_k(1)=0,\ \ \ff k\ne0.
\earr\end{equation}

Now, let us check \eqref{e4.9} or, equivalently, \eqref{e4.10}. We have for $u=u^j$
$$\derpp n\ u(x,y)=-i\dd\sum_k\frac{e^{ikx}}k v''_k(y),\ \ \ff x,y\in0,1,$$and so \eqref{e4.10} reduces to
\begin{equation}\label{e4.16}
|v''_k(0)|+|v''_k(1)|>0\mbox{\ \ for all $k$.}\end{equation}Assume that $v''_k(0)=v''_k(1)=0$ for all $k$ and lead from this to a contradiction.  To this end we set $W_k=v''_k-k^2v_k$ and rewrite \eqref{e4.12} as
\begin{equation}\label{e4.15nou}
\barr{c}
-\nu W''_k+(\nu k^2+ikU-\lbb)W_k=ikU''v_k\ \inr\ (0,1),\vsp
W_k(0)=W_k(1)=0.\earr \end{equation}

If we multiply \eqref{e4.15nou} by $\ov W_k$,  integrate on $(0,1)$ and take the real part, we obtain that
$$\int^1_0(\nu|W'_k|^2+(\nu k^2-\Res\,\lbb)|W_k|^2)dy=0,\ \ff k$$
and since $\Res\,\lbb=\Res\,\lbb_j\le0$ for all $j=1,...,N$, we get $W_k\equiv0$,   and so $v_k\equiv0.$
 The contradiction we arrived at proves \eqref{e4.16}
 and \eqref{e4.8}. Arguing similarly for the dual system with eigenfunctions $\{u^*_k,v^*_k\}$, we get
 for $W^*_k=(v^*_k)''-k^2v^*_k$
$$\barr{c}
L^*_kv^*_k\equiv-\nu(W^*_k)''+(\nu k^2-ikU-\ov\lbb)W^*_k-2ik(v^*_k)'U'=0\ \mbox{ in }(0,1),\vsp
  v^*_k(0)=v^*_k(1)=0,\ \ (v^*_k)'(0)=(v^*_k)'(1)=0,\ \ (v^*_k)''(0)=(v^*_k)''(1)=0,\earr$$
and argue from  this to a contradiction. We set $\calx^*=\{\vf;\ L^*_k\vf=0,\ \vf(0)=\vf(1)=0,\ \vf'(0)=\vf'(1)=0\}$  and let $\check\vf(y)=\vf(1-y),$ $y\in[0,1]$. Since $\dim\calx^*\le2$ and $\check U\equiv U$, we infer that each $\vf\in\calx^*$ is either symmetric (that is, $\vf\equiv\check\vf)$ or antisymmetric (that is, $\vf\equiv-\check\vf)$. Assume that $v^*_k$ is symmetric. By \eqref{e4.15nou} we see via integration by parts that $\int^1_0|v^*_k|^2dy=0$ if there is $\vf$ such that
\begin{equation}\label{e4.15aa} L_k\vf=v^*_k\ \mbox{ in }(0,1);\ \ \ \vf(0)=\vf(1)=0.\end{equation}
In order to show that there is such a function $\vf$, we shall prove that there is $\vf_1$ such that
\begin{equation}\label{e4.15aaa} L_k\vf_1=0\mbox{ in }(0,1),\ \ \ \ \vf_1(0)+\vf_1(1)\ne0.\end{equation}
Indeed, if such a $\vf_1$ exists, by replacing $\vf_1$ by $\vf+\check\vf_1$, we may assume that $\vf_1$ is symmetric. If $\vf_2$ is a symmetric solution to $L_k\vf_2=v^*_k$, then clearly $\vf=\vf_2-\vf_2(0)(\vf_1(0))\1\vf_1$ satisfies \eqref{e4.15aa} because, by \eqref{e4.15aaa}, $\vf_1(0)\ne0$.

Now, to prove the existence in \eqref{e4.15aaa}, we shall argue as in \cite{18} and assume that $\calx=\{\psi;\ L_k\psi=0\}\equiv\{\psi;\ L_k\psi=0,\ \psi(0)+\psi(1)=0\}$ and argue from this to a contradiction. We set $\calx_1=\{\psi\in\calx;\ \psi''(0)+\psi''(1)=0\}$ and prove that $\psi=-\check\psi$ for each $\psi\in\calx_1$. Indeed, $\theta=\psi+\check\psi$ satisfies $\theta(0)=\theta(1)=0$, $\theta''(0)=\theta''(1)=0$ and $W=\theta''-k^2\theta$ satisfies \eqref{e4.15nou} with $v_k=\theta$. Then, we obtain as above  that $W\equiv0,$ $\theta\equiv0.$ The spaces $\calx_1=\{\psi\in\calx;\ \psi\equiv-\check\psi\}$ and $\calx_2=\{\psi\in\calx;\ \psi\equiv\check\psi\}$           are orthogonal and both have  dimension 2 because
$$\barr{lcl}
\calx_1&=&\{\psi\in\calx;\ \psi'(\frac12)=\psi'''(\frac12)=0\},\vsp
\calx_2&=&\{\psi\in\calx;\ \psi(\frac12)=\psi''(\frac12)=0\}.\earr$$Hence, $\calx=\calx_1\oplus\calx_2$. On the other hand, the space $\{\psi\in\calx;\ \psi''(0)=0\}$ which has dimension 3, has nonempty intersection with $\calx_2$. Hence, there is $\psi\in\calx$ symmetric such that $\psi''(0)=\psi''(1)=0$. Clearly, $\psi\in\calx_1$, which is absurd. This completes the proof.
\bigskip

As regards (H1), it is not clear if it is always satisfied in the present situation and so we keep it.

\begin{lemma}\label{l2} The systems
$\left\{\derp{\vf_j}n\right\}^N_{j=1}$ and
$\left\{\derp{\vf^*_j}n\right\}^N_{j=1}$  are linearly independent
on $\pp\calo$.\end{lemma}

\pf It suffices to prove the independence of $\lac\derp{\vf_j}n\rac^N_{j=1}$, $\vf_j=(u^j,v^j)$. We have as above $u^j=\{u^j_k\}_k,$ $v^j=\{v^j_k\}$, $j=1,...,N.$ If $\{\vf_j\}$ are eigenvectors corresponding to the same eigenvalue, the independence follows by Lemma \ref{l1}. Assume that $\cala\vf_1=\lbb_1\vf_1$, $\cala\vf_2=\lbb_2\vf_2$, where $\lbb_1\ne\lbb_2$, and that
$$\derp{\vf_1}n+  \derp{\vf_2}n=0\ \inr\ y=0,1.$$
We set $\wt v_k=v^1_k+  v^2_k$ and $\wt W_k=\wt v''_k-k^2\wt v_k.$ Then, we have as above (see \eqref{e4.15nou})
$$\barr{l}
-\nu\wt W''_k+(\nu k^2+ikU-\lbb_1)\wt W_k-(\lbb_2-\lbb_1)((v^2_k)''-k^2v^2_k)
=ikU''\wt v_k,\vsp
-\nu\wt W''_k+(\nu k^2+ikU-\lbb_2)\wt W_k-(\lbb_1-\lbb_2)((v^1_k)''-k^2v^1_k)=ikU''\wt v_k.
\earr$$
This yields
$$\barr{l}
\dd\int^1_0(\nu|\wt W'_k|^2+(\nu k^2-\Res\,\lbb_1)|\wt W_k|^2)dy-\Res
\left[(\lbb_2-\lbb_1)\dd\int^1_0((v^2_k)''-k^2v^2_k)\ov{\wt W}_k)dy\right]=0,\vsp
\dd\int^1_0(\nu|\wt W'_k|^2+(\nu k^2-\Res\,\lbb_2)|\wt W_k|^2)dy-
\Res\left[(\lbb_1-\lbb_2) \dd\int^1_0(v^1_k)''-k^2v^2_k)\wt W_kdy\right]=0,\earr$$
 and, therefore,
 $$\int^1_0(\nu|\wt W'_k|^2+(\nu k^2-\Res\,\lbb_1-\Res\,\lbb_2)|\wt W_k|^2)dy=0.$$
 Hence, $\wt W_k\equiv0$, $\wt v_k=0$, which is absurd because $v^1_k,v^2_k$ are independent.

 By induction with respect to $j$, one proves the independence of $\lac\derp{\vf_j}n\rac^N_{j=1}$.  The case $\lac\derp{\vf^*_j}n\rac^N_{j=1}$ is completely similar.

Now, following the general case \eqref{e3.5}, we can design a feedback controller $(u^0,v^0)$ for system \eqref{e4.3}, \eqref{e4.5}. We set
$$\vf^*_j=(u^*_j,v^*_j),\ \ j=1,...,N,$$where $\vf^*_j$ are eigenvectors of the dual operator $\cala^*$ with corresponding eigenvalues $\ov\lbb_j$ and $\Res\,\lbb_j<0$ for $j=1,...,N.$

We consider the feedback controller
\begin{equation}\label{e4.19}
\barr{rcll}
u^0(t,x,y)&{=}&\eta\dd\sum^N_{j=1}\mu_j v_j(t)\phi^1_j(x,y),\ x\in\rr,\ y=0,1,\vsp
v^0(t,x,y)&{=}&\eta\dd\sum^N_{j=1}\mu_j v_j(t)(\phi^2_j(x,y)+\a H(y)),\ x\in\rr,\ y=0,1,\vsp
v_j(t)&{=}&\dd\int^{2\pi}_0
(u(t,x,y)u^*_j(x,y)+v(t,x,y)v^*_j(x,y))dx\,dy\vsp
&{=}&\dd\sum_k(u_k(t,y)(\bar u^*_j)_k(y)+v_k(t,y)(\bar v^*_k)_k(y)).\earr\end{equation}Here $\a$ is an arbitrary constant, $H(0)=-1$, $H(1)=1$, $\mu_j$ are defined as \eqref{e2.6} and, according to \eqref{e2.7}, $\phi^i_j$, $i=1,2,$ are of the form
$$\phi^1_j=\sum^N_{j=1}\a_{ij}(u^*_i)'(y),\ \ \phi^2_j=\dd\sum^N_{i=1}\a_{ij}(v^*_j)'(y),$$where $\a_{ij}$ are chosen as in Section 2.2. Then, by Theorem \ref{t1}, we have

\begin{theorem}\label{t3a} For each $\a\in\rr$ and $\eta$ suitable
chosen, the feedback boundary controller \eqref{e4.19} stabilizes
exponentially system \eqref{e4.3}.\end{theorem}

We note that condition \eqref{e1.4} automatically holds in this
case for any constant $\a$. However, by Theorem \ref{t2}, it
follows also the stabilization with a controller $(u^0,v^0)$ with
support in $\{y=0\}$ or $\{y=1\}$ if $\a=\a(x,y)$ is taken in such
a way that $\dd\int^{2\pi}_0\a(x)dx=0$.

We note also that, by Theorem \ref{t3a}, we infer that the feedback
controller \eqref{e4.19} is exponentially stabilizable in the \NS\
equation \eqref{e4.1}.

\begin{remark}\label{r4.1} {\rm The boundary stabilization of \eqref{e4.1} was studied in \cite{1}, \cite{2}, \cite{3}, \cite{5}, \cite{20}, \cite{21}. In \cite{3} and \cite{16a} it is proved the existence of a normal stabilizing controller $\{u_k,v_k\}$ such that $u_k\equiv v_k\equiv0$ for $|k|\ge M,$ which is, apparently, a stronger result than Theorem \ref{t3a}. However, the advantage of the present result is the explicit design of the feedback controller.

Note also that, by Theorem \ref{t3}, the feedback controller is stabilizable in Navier--Stokes equation \eqref{e4.1}. Also, as in Theorem \ref{t4}, it can be replaced by a real feedback controller.}\end{remark}

\end{document}